\documentstyle{amsppt}
\magnification=1100
\NoBlackBoxes
\advance\hsize-15pt

\topmatter
\nopagenumbers
\font\titlefont=cmbx12\font\authorfont=cmr10
\baselineskip=13pt
\rightheadtext{MOORE--PENROSE INVERSE}
\leftheadtext{E. TEVELEV}
\topmatter
\centerline{\titlefont
MOORE--PENROSE INVERSE,}
\centerline{\titlefont
PARABOLIC SUBGROUPS,}
\centerline{\titlefont
AND JORDAN PAIRS}
\vskip 1.2pc
\centerline{\authorfont
E. TEVELEV$^{\ast}$}
\vskip 1.1pc
\centerline{Moscow Independent University}
\smallskip
\centerline{tevelev\@mccme.ru}

\thanks{$^\ast$ 
The research was supported by the grant
INTAS-OPEN-97-1570 of the INTAS foundation.}\endthanks

\abstract
A Moore--Penrose inverse of an arbitrary complex matrix $A$
is defined as a unique matrix $A^+$ such that
$AA^+A=A$, $A^+AA^+=A^+$, and $AA^+$, $A^+A$ are Hermite matrices.
We show that this definition has a natural generalization
in the context of shortly graded simple Lie algebras corresponding
to parabolic subgroups with {\it aura} (abelian unipotent radical)
in simple complex Lie groups,
or equivalently in the context of simple complex Jordan pairs.
We give further generalizations and applications.
\endabstract

\endtopmatter
\def\cal#1{{\fam2 #1}}
\def\C{{\Bbb C}}
\def\H{{\Bbb H}}
\def\R{{\Bbb R}}
\def\Ca{\Bbb C\text{a}}
\def\OO{{\Bbb O}}
\def\A{{\Bbb A}}
\def\g{\goth{g}}
\def\h{\goth{h}}
\def\z{\goth{z}}
\def\ssl{\goth{sl}}
\def\sp{\goth{sp}}
\def\su{\goth{su}}
\def\so{\goth{so}}

\def\k{\goth{k}}
\def\p{\goth{p}}
\def\q{\goth{q}}
\def\l{\goth{l}}
\def\n{\goth{n}}
\def\t{\goth{t}}
\def\p{\goth{p}}
\def\u{\goth{u}}
\def\End{{\text{\rm End}}}
\def\Ker{{\text{\rm Ker}}}

\def\Im{{\text{\rm Im}}}

\def\Aut{{\text{\rm Aut}}}

\def\GL{\text{GL}}
\def\Tr{\text{Tr}}
\def\SL{\text{SL}}

\def\SO{\text{SO}}

\def\Sp{\text{Sp}}
\def\ad{\text{ad}}
\def\Ad{\text{Ad}}
\def\Ann{\text{Ann}}
\def\Mat{\text{Mat}}
\def\Z{\Bbb Z}
\def\O{\cal O}
\def\ow{\goth{o}(\omega)}
\def\Ow{\text{O}(\omega)}
\def\hht{\text{ht}}
\def\rk{\text{rk}}
\def\Hom{\text{Hom}}
\def\Id{\text{Id}}

\document
\head Introduction\endhead
The nice notion of a generalized inverse of an arbitrary matrix
(possibly singular or even non-square) has been discovered
independently by Moore \cite{Mo} and Penrose \cite{Pe}.
The following 
definition belongs to Penrose (Moore's definition is different but 
equivalent):

\definition{Definition}
A matrix $A^+$ is called {\it a MP-inverse\/} of a 
matrix $A$
if 
$$AA^+A=A,\quad A^+AA^+=A^+,$$ 
and $AA^+$, $A^+A$ are Hermite matrices.
\enddefinition

It is quite surprising but a MP-inverse always exists and is unique.
Since the definition is symmetric with respect to $A$ and $A^+$
it follows that $(A^+)^+=A$.
If $A$ is a non-singular square matrix then $A^+$
coincides with an ordinary inverse matrix $A^{-1}$.
The theory of MP-inverses and their numerous modifications becomes
now a separate subfield of Linear Algebra \cite{CM} with various applications.
The aim of this paper is to demonstrate that this notion
quite naturally arises in the theory of shortly graded simple Lie algebras.
To explain this connection let us first give another definition of a MP-inverse.

\definition{Equivalent definition of a MP-inverse}
Suppose that $A\in\Mat_{n,m}(\C)$. 
Then a matrix $A^+\in\Mat_{m,n}(\C)$
is called {\it a MP-inverse\/} of $A$ if there exist Hermite matrices 
$B_1\in\Mat_{n,n}(\C)$ and $B_2\in\Mat_{m,m}(\C)$
such that the following matrices form an $\ssl_2$-triple
in $\ssl_{n+m}(\C)$:
$$E=\left(\matrix 0&A\cr0&0\cr\endmatrix\right),\quad 
H=\left(\matrix B_1&0\cr0&B_2\cr\endmatrix\right),\quad
F=\left(\matrix 0&0\cr A^+&0\cr\endmatrix\right).$$
\enddefinition

\remark{Remark}
By an $\ssl_2$-triple $\langle e,h,f\rangle$ in a Lie algebra $\g$
we mean a collection of (possibly zero) vectors
such that 
$$[e,f]=h,\quad [h,e]=2e,\quad [h,f]=-2f.$$
In other words, an $\ssl_2$-triple is a homomorphic image
of canonical generators of $\ssl_2$ with respect to some homomorphism 
of Lie algebras $\ssl_2\to\g$.
\endremark

This definition admits an immediate generalization.
In the sequel we shall use various facts about shortly graded
simple Lie algebras without specific references to original papers,
the reader may consult, for example, 
papers \cite{RRS}, \cite{Pa}, or \cite{MRS}
for explanations and further references.
All necessary facts about complex and real Lie groups, 
Lie algebras, and algebraic groups
can be found in \cite{VO}.

Suppose that $\g$ is a simple complex Lie algebra, $G$
is a corresponding simple simply-connected Lie group.
Suppose further that $P$ is a parabolic subgroup of $G$
with abelian unipotent radical (with {\it aura}).
Then $\g$ admits a short grading
$$\g=\g_{-1}\oplus \g_0\oplus\g_1$$
with only three nonzero parts. Here $\p=\g_0\oplus\g_1$ is a Lie
algebra of $P$ and $\exp\g_1$ is the abelian unipotent radical of $P$.
Let $\k_0$ be a compact real form of $\g_0$.

\remark{Remark}
In this paper
we shall permanently consider compact real forms of reductive
subalgebras of simple Lie algebras. These subalgebras
will always be Lie algebras of algebraic reductive subgroups
of a corresponding simple complex algebraic group. 
Their compact real forms will always be understood
as Lie algebras of compact real forms of corresponding
algebraic groups. 
For example, a Lie algebra of an algebraic torus has a unique
compact real form.
\endremark

Suppose now that $e\in\g_1$. It is well-known that
there exists a {\it homogeneous} $\ssl_2$-triple $\langle e,h,f\rangle$
such that $h\in\g_0$ and $f\in\g_{-1}$.

\definition{Definition}
An element $f\in\g_{-1}$ is called {\it a MP-inverse\/} of $e\in\g_1$
if there exists a homogeneous $\ssl_2$-triple $\langle e,h,f\rangle$
with $h\in i\k_0$.
\enddefinition

MP-inverses of elements $f\in\g_{-1}$ are defined in the same
way. It is clear that if $f$ is a MP-inverse of $e$
then $e$ is a MP-inverse of $f$. 

\example{Example}
Suppose that $G=\SL_{n+m}$ and $P\subset G$ is a maximal
parabolic subgroup of block triangular matrices
of the form 
$$\left(\matrix B_1&A\cr 0&B_2\cr\endmatrix\right),\quad\text{where}\quad
B_1\in\Mat_{n,n},\ A\in\Mat_{n,m},\ B_2\in\Mat_{m,m}.$$
The graded components of the correspondent grading
consist of matrices of the following form:
$$\g_{-1}=\left(\matrix 0&0\cr A'&0\cr\endmatrix\right),\quad
\g_0=\left(\matrix B_1&0\cr 0&B_2\cr\endmatrix\right),\quad
\g_1=\left(\matrix 0&A\cr 0&0\cr\endmatrix\right),$$
where $A'\in\Mat_{m,n}$, $B_1\in\Mat_{n,n}$, $B_2\in\Mat_{m,m}$, and
$A\in\Mat_{n,m}$.
One can take $\k_0$ to be a real Lie algebra of block diagonal
skew-Hermite matrices with zero trace. Then $i\k_0$ is a vector space
of block diagonal Hermite matrices with zero trace.
Therefore in this case we return to a previous definition
of a Moore--Penrose inverse.
\endexample

Our first result is the following
\proclaim{Theorem 1}
For any $e\in\g_1$ there exists a unique MP-inverse $f\in\g_{-1}$.
\endproclaim

It obviously follows that
for any non-zero $f\in\g_{-1}$ there exists a unique MP-inverse $e\in\g_1$.
So taking a MP-inverse is a well-defined involutive operation.
In general, it is not equivariant with respect to a Levi subgroup $L\subset P$
with Lie algebra $\g_0$,
but only with respect to its maximal compact subgroup $K_0\subset L$.

Theorem~1 will be proved in \S2 by a general argument,
without using case-by-case considerations. But the 
classification of parabolic subgroups with aura
in simple groups is, of course, well-known.
We have tried to give an intrinsic description
of the Moore-Penrose inverse in all arising cases.
The calculation of Moore--Penrose inverses
arising from short gradings of classical simple Lie algebras
is quite straight-forward, so
we shall give here only the summary of these calculations and
avoid proofs.

\subhead Linear maps\endsubhead
This is, of course, the classical Moore--Penrose inverse.
Let us recall its intrinsic description.
Suppose that $\C^n$ and $\C^m$ are vector spaces
equipped with standard Hermite scalar products.
For any linear map $F:\,\C^n\to\C^m$ its Moore-Penrose
inverse is a linear map $F^+:\,\C^m\to\C^n$ defined as follows.
Let $\Ker F\subset\C^n$ and $\Im F\subset \C^m$ be
a kernel and an image of $F$.
Let  $\Ker^\perp F\subset\C^n$ and $\Im^\perp F\subset \C^m$
be their orthogonal complements with respect 
to the Hermite scalar products.
Then $F$ defines via restriction a bijective linear map
$\tilde F:\,\Ker^\perp F\to\Im F$.
Then $F^+:\,\C^m\to\C^n$ is a unique linear map such that
$F^+|_{\Im^\perp F}=0$ and 
$F^+|_{\Im F}=\tilde F^{-1}$.
This MP-inverse corresponds to short gradings of $\ssl_{n+m}$.

\subhead Symmetric and skew-symmetric bilinear forms\endsubhead
Suppose that $V=\C^n$ is a vector space equipped with 
a standard Hermite scalar product.
For any symmetric (resp.~skew-symmetric) bilinear form $\omega$ on $V$
its Moore-Penrose inverse is 
a symmetric (resp.~skew-symmetric) bilinear form $\omega^+$ on $V^*$
defined as follows.
Let $\Ker\, \omega\subset V$ be the kernel of~$\omega$.
Then $\omega$ induces a non-degenerate bilinear form
$\tilde\omega$ on $V/\Ker\,\omega$.
Let $\Ann(\Ker\, \omega)\subset V^*$ 
be an annihilator of $\Ker\, \omega$. 
Then $\Ann (\Ker\, \omega)$ is canonically isomorphic to the dual of $V/\Ker\,\omega$.
Therefore the form
$\tilde\omega^{-1}$ on 
$\Ann (\Ker\, \omega)$ is well-defined.
The form $\omega^+$ is defined as a unique form such that
its restriction on
on $\Ann(\Ker\,\omega)$ coincides with 
$\tilde\omega^{-1}$ and its kernel is $\Ann(\Ker\,\omega)^\perp$,
the orthogonal complement with respect to a standard
Hermite scalar product on $V^*$.
This MP-inverse corresponds to the short grading of $\sp_{2n+2}$ 
(resp.~$\so_{2n+2}$).

\subhead Vectors in a vector space with scalar product\endsubhead
Let $V=\C^n$ be a vector space with standard bilinear
scalar product $(\cdot,\cdot)$.
For any vector $v\in V$ its Moore-Penrose inverse
$v^\vee$ is again a vector in $V$ defined as follows:
$$v^\vee=\cases
{\textstyle2v\over\textstyle (v,v)},&\text{if}\ (v,v)\ne0\cr
{\textstyle\overline v\over\textstyle (\overline v,v)},&\text{if}\ (v,v)=0,\ v\ne0\cr
0,&\text{if}\ v=0.\cr
\endcases$$
This MP-inverse corresponds to the short grading of~$\so_{n+2}$.

The short gradings of exceptional Lie algebras $E_6$ and $E_7$
deserve more detailed considerations.
This is done in \S3.
It is well-known that the theory of shortly graded simple Lie algebras
is equivalent to the theory of finite-dimensional simple Jordan pairs.
It turns out that the Moore--Penrose inversion has a very
simple interpretation in this alternative language.
We describe this connection also in \S3.

In \S6 we consider shortly graded real simple Lie algebras.
It turns out that the analogue of Theorem~1 is also
true in this case. 

It is quite natural to ask whether it is possible
to extend the notion of the Moore--Penrose inverse
from parabolic subgroups with aura to arbitrary parabolic subgroups.
It is also interesting to consider the ``non-graded'' situation.
Let us start with it. Suppose $G$ is a simple connected simply-connected
Lie group with Lie algebra $\g$. We fix a compact real form $\k\subset\g$.

\definition{Definition}
A nilpotent orbit $\O\subset\g$ is called {\it a Moore--Penrose orbit\/} if for any
$e\in\O$ there exists an $\ssl_2$-triple $\langle e,h,f\rangle$
such that $h\in i\k$.
\enddefinition

It turns out that it is quite easy to find all Moore-Penrose orbits.
Recall that the height $\hht(\O)$ of a nilpotent orbit $\O=\Ad(G)e$
is equal to the maximal integer $k$ such that $\ad(e)^k\ne0$.
Clearly $\hht(\O)\ge2$.

\proclaim{Theorem 2}
$\O$ is a Moore--Penrose orbit if and only if $\hht(\O)=2$.
In this case for any $e\in\O$ there exists a unique $\ssl_2$-triple
$\langle e,h,f\rangle$ such that $h\in i\k$.
\endproclaim

This theorem will be proved in \S1. It is worthy to mention
here the following result of Panyushev \cite{Pa1}:
$\hht(\O)\le3$ if and only if $\O$ is a spherical $G$-variety
(that is, a Borel subgroup $B\subset G$ has an open orbit in $\O$).
Therefore, all Moore--Penrose orbits are spherical.
If $G=\SL_n$ or $G=\Sp_n$ then the converse is also true.

Now let us turn to the graded situation.
Suppose that $\g$ is a $\Z$-graded simple Lie algebra, 
$\g=\mathop{\oplus}\limits_{k\in\Z}\g_k$.
Let $P\subset G$ be a parabolic subgroup with Lie algebra 
$\p=\mathop{\oplus}\limits_{k\ge0}\g_k$. Let $L\subset P$ be a Levi subgroup with Lie algebra
$\g_0$. 
We choose a compact real form $\k_0$ of $\g_0$.
Suppose now that $e\in\g_k$. It is well-known that
there exists a homogeneous $\ssl_2$-triple $\langle e,h,f\rangle$
with $h\in\g_0$ and $f\in\g_{-k}$.

\definition{Definition}
Take any $k>0$ and any $L$-orbit $\O\subset\g_k$.
Then $\O$ is called {\it a Moore-Penrose orbit\/} if for any $e\in\O$
there exists a homogeneous $\ssl_2$-triple $\langle e,h,f\rangle$
such that $h\in i\k_0$. In this case $f$ is called a MP-inverse of $e$.
A grading is called {\it
a Moore--Penrose grading
in degree $k>0$\/} if all $L$-orbits in $\g_k$ are Moore-Penrose. 
A grading is called {\it a Moore--Penrose grading\/} if it
is a Moore--Penrose grading in any positive degree.
A parabolic subgroup $P\subset G$ is called {\it a Moore--Penrose
parabolic subgroup\/} if there exists a Moore--Penrose grading 
$\g=\mathop{\oplus}\limits_{k\in\Z}\g_k$ such that 
$\p=\mathop{\oplus}\limits_{k\ge0}\g_k$ is a Lie algebra of $P$.
\enddefinition

One should be careful comparing graded and non-graded situation:
if $\O\subset\g_k$ is a Moore--Penrose $L$-orbit then
$\Ad(G)\O\subset\g$ is not necessarily a Moore--Penrose
$G$-orbit.
Let us give a criterion 
for an $L$-orbit to be Moore--Penrose.
Suppose $\O=\Ad(L)e\subset\g_k$. Take any homogeneous $\ssl_2$-triple
$\langle e,h,f\rangle$. Then $h$ defines a grading 
$\g_0=\mathop{\oplus}\limits_{n\in\Z}\g_0^n$,
such that $\ad(h)|_{\g_0^n}=n\cdot\Id$.

\proclaim{Theorem 3}
$\O$ is a Moore--Penrose orbit if and only if 
$\ad(e)g_0^n=0$ for any $n>0$.
In this case for any $e'\in\O$ there exists a unique homogeneous 
$\ssl_2$-triple $\langle e',h',f'\rangle$ such that
$h'\in i\k_0$.
\endproclaim

This Theorem will be proved in \S2. 
It gives a characterization of Moore--Penrose orbits
independent on the choice of a compact form and also
provides an algorithm for checking the Moore--Penrose property.

It is easy to see that in the graded
situation a Moore--Penrose orbit is not necessarily spherical.
However, some interesting orbits are both spherical and Moore-Penrose.
Let us give several examples.

\example{Example 1}
If $P$ is a parabolic subgroup with aura then all $L$-orbits in $\g_1$
are Moore--Penrose by Theorem~1. It is well-known that all of them are also
spherical. More generally, take any grading of $\g$ and
suppose that $d$ is equal 
to the maximal $k$ such that $\g_k\ne0$ (the {\it height\/} of grading).
Then all $L$-orbits in $\g_k$
are both spherical and Moore-Penrose for $k>d/2$. This fact easily follows
from the previous remark. (Consider the short-graded Lie algebra
$\g_{-k}\oplus\g_0\oplus\g_k$. Of course it is not necessarily simple
but this is not essential.)
\endexample

\example{Example 2}
Suppose that $G$ is a simple group of type $\text{G}_2$.
We fix a root decomposition.
There are two simple roots $\alpha_1$ and $\alpha_2$ such that
$\alpha_1$ is short and $\alpha_2$ is long.
There are $3$ proper parabolic subgroups: Borel subgroup $B$
and two maximal parabolic subgroups $P_1$ and $P_2$ such that
a root vector of $\alpha_i$ belongs to a Levi subgroup of $P_i$.
Then the following is an easy application of Theorem~3.
$B$ is a Moore--Penrose parabolic subgroup (actually
Borel subgroups in all simple groups are Moore--Penrose parabolic
subgroups
with respect to any grading).
$P_1$ is not Moore--Penrose, but it is a Moore--Penrose
parabolic subgroup in degree $2$ (with respect to the
natural grading of height~$2$). 
$P_2$ is a Moore--Penrose parabolic subgroup.
\endexample

\example{Example 3}
Suppose $G=\SL_n$.
We fix positive integers $d_1,\ldots,d_k$ such that $n=d_1+\ldots+d_k$.
We consider the parabolic subgroup $P(d_1,\ldots,d_k)\subset\SL_n$
that consists of all upper-triangular block matrices with sizes
of blocks equal to $d_1,\ldots,d_k$. We take a standard grading.
Then $\g_1$ is identified with the linear space of all tuples of linear maps
$\{f_1,\ldots,f_k\}$,
$$\C^{d_1}\mathop{\longleftarrow}\limits^{f_1}\C^{d_2}
\mathop{\longleftarrow}\limits^{f_2}\ldots
\mathop{\longleftarrow}\limits^{f_{k-1}}\C^{d_k},$$
$\g_{-1}$ is identified with the linear space of all tuples of linear maps
$\{g_1,\ldots,g_k\}$,
$$\C^{d_1}\mathop{\longrightarrow}\limits^{g_1}\C^{d_2}
\mathop{\longrightarrow}\limits^{g_2}\ldots
\mathop{\longrightarrow}\limits^{g_{k-1}}\C^{d_k},$$
and Levi subgroup $L(d_1,\ldots,d_k)$
is just a group of all $k$-tuples 
$$(A_1,\ldots,A_k)\in\GL_{d_1}\times\ldots
\times\GL_{d_k}\quad\text{such that}\quad\det(A_1)\cdot\ldots\cdot\det(A_k)=1,$$
acting on these spaces of linear maps in an obvious way.
The most important among $L$-orbits 
are varieties of complexes.
To define them, let us fix in addition non-negative integers
$m_1,\ldots,m_{k-1}$ such that $m_{i-1}+m_i\le d_i$ (we set $m_0=m_k=0$),
and consider the subvariety of all tuples $\{f_1\ldots,f_{k-1}\}$
as above such that $\rk f_i=m_i$ and $f_{i-1}\circ f_i=0$ for any $i$.
These tuples form a single $L$-orbit~$\O$ called a variety of complexes.
It is well-known that $\O$ is spherical.
For any tuple $\{f_1,\ldots,f_{k-1}\}\in\O$
consider the tuple $\{f_1^+,\ldots,f_{k-1}^+\}\in\g_{-1}$,
where $f_i^+$ is a classical ``matrix'' Moore--Penrose inverse of $f_i$.
The reader may check that this new tuple is again a complex,
moreover, this complex is a Moore--Penrose inverse
(in our latest meaning of this word)
of an original complex.
In particular, orbits of complexes are Moore--Penrose orbits.
\endexample

From the first glance only few parabolic subgroups are Moore--Penrose.
But this is scarcely true.
For example, we have the following Theorem:

\proclaim{Theorem 4}
Any parabolic subgroup in $SL_n$ is Moore--Penrose.
\endproclaim

This Theorem will be proved in \S4. We shall also describe
an algorithm there, which shows that in order to find
all Moore--Penrose parabolic subgroups in some simple group $G$
it is sufficient to determine all Moore--Penrose maximal parabolic 
subgroups in simple components of Levi subgroups of $G$.
In particular in order to find all Moore--Penrose parabolic
subgroups in classical simple groups it suffices
to do this job only for maximal parabolic subgroups.
We shall do this also in \S4.

To explain our interest in Moore--Penrose parabolic subgroups let us reproduce
a conjecture from \cite{Te}.
Suppose once again that $G$ is a simple connected simply-connected
Lie group, $P$ is its parabolic subgroup, $\p\subset\g$
are their Lie algebras.
We take any irreducible $G$-module $V$.
There exists a unique proper $P$-submodule $M_V$ of $V$.
We have the inclusion $i:\,M_V\to V$, the projection
$\pi:\,V\to V/M_V$ and the map $R_V:\,\g\to\End(V)$ defining the representation.
Therefore we have a linear map $\tilde R_V:\,\g\to\Hom(M_V,V/M_V)$, namely
$\tilde R_V(x)=\pi\circ R_V(x)\circ i$. Clearly $\p\subset\Ker\tilde R_V$.
Therefore, we finally have a linear map $\Psi_V:\,\g/\p\to\Hom(M_V,V/M_V)$.

\proclaim{Conjecture}
There exists an algebraic stratification $\g/\p=\mathop{\sqcup}\limits_{i=1}^nX_i$
such that for any $V$ the function $\rk\,\Psi_V$ is constant
along each $X_i$.
\endproclaim

These stratifications were used in \cite{Te} in order
to solve some geometric problems similar to the classical
problem of determining the maximal dimension of a projective
subspace contained in a generic hypersurface of a given degree
in a projective space.

It is clear that all functions 
$\rk\,\Psi_V$ are $P$-invariant. Therefore if $P$ has finitely many
orbits in $\g/\p$ then the conjecture is true.
By Pyasetsky theorem \cite{P} this holds if and only if 
$P$ has finitely many orbits in the dual module $(\g/\p)^*$,
or, equivalently, in the unipotent radical of $P$.
All parabolic subgroups with this property are now completely classified \cite{HR}.
There are not too many of them.
It turns out that there is another case when the conjecture
is true.

\proclaim{Theorem 5}
Suppose that 
a grading $\g=\mathop{\oplus}\limits_{k\in\Z}\g_k$ is a Moore--Penrose grading
in all positive degrees except at most one.
Then the Conjecture is true for the corresponding parabolic subgroup $P$.
\endproclaim

This Theorem is proved in \S5. 
For example, combining Theorem~4, Theorem~5, and Example~2
we get the following corollary:

\proclaim{Corollary}
The conjecture is true for any parabolic subgroup in $\SL_n$ or $\text{G}_2$.
\endproclaim

This paper was written during my stay in the Mathematical Institute
in Basel.
I would like to thank prof.~H.~Kraft for the warm hospitality.

\head \S1. Moore--Penrose orbits in simple Lie algebras\endhead
In this section we prove Theorem~2. Recall that
$G$ is a simple connected simply-connected Lie group with a Lie algebra $\g$
and a compact form $\k\subset\g$.
Let $x\to\overline x$ denotes a complex conjugation in $\g$ with respect to 
the compact form $\k$. 
Therefore $x=\overline x$ iff $x\in\k$ and $x=-\overline x$ iff $x\in i\k$.
Let $B(x,y)=\Tr\, \ad(x)\ad(y)$ be the Killing form of $\g$.
Finally, 
let $H(x,y)=-B(x,\overline y)$ be a positive-definite Hermite form on $\g$.

\proclaim{Lemma 1.1}
We fix a nilpotent element $e\in\g$.
Suppose that $\langle e,h,f\rangle$ is an $\ssl_2$-triple in $\g$ such that 
$h\in i\k$. Then for any other $\ssl_2$-triple 
$\langle e,h',f'\rangle$ we have $H(h,h)<H(h',h')$.
In particular,
if there exists an $\ssl_2$-triple $\langle e,h,f\rangle$
with $h\in i\k$ then the $\ssl_2$-triple with this property is unique.
\endproclaim 

\demo{Proof}
Recall that if $\langle e,h,f\rangle$ is an $\ssl_2$-triple then 
$h$ is called {\it a characteristic\/} of~$e$.
Consider the subset $\cal H\subset\g$ consisting 
of all possible characteristics of~$e$. It is well-known that
$\cal H$ is an affine subspace in $\g$ such that the corresponding
linear subspace is precisely the unipotent radical $\z^u_{\g}(e)$ of the centralizer
$\z_{\g}(e)$ in~$\g$ of the element $e$.
Since $H(h',h')$ is a strongly convex function on $\cal H$, 
there exists a unique element $h_0\in\cal H$
such that $H(h_0,h_0)<H(h',h')$ for any $h'\in\cal H$, $h'\ne h_0$.
We need to show that $h_0=h$.
It is clear that an element $h_0\in\cal H$ minimizes 
$H(h,h)$ on $\cal H$ iff $H(h_0,\z^u_{\g}(e))=0$ iff
$B(\overline h_0,\z^u_{\g}(e))=0$.
If $h\in \cal H\cap i\k_0$ then $\overline h=-h$ and we have
$$B(\overline h,\z^u_{\g}(e))=
-B(h,\z^u_{\g}(e))=
-B([e,f],\z^u_{\g}(e))=
B(f,[e,\z^u_{\g}(e)])=0.$$
Therefore $h=h_0$.
\qed\enddemo

Suppose that $\langle e,h,f\rangle$ is an $\ssl_2$-triple in $\g$.
Consider the grading $\g=\mathop{\oplus}\limits_k\g_k$ such that $x\in\g_k$
iff $[h,x]=kx$. Let $\n_+=\mathop{\oplus}\limits_{k>0}\g_k$, $\n_-=\mathop{\oplus}\limits_{k<0}\g_k$.
It is well known that $\z_\g^u(e)\subset\n_+$.

\proclaim{Lemma 1.2}
Suppose that $\z_\g^u(e)=\n_+$. Then $\O=\Ad(G)e$ is a Moore--Penrose orbit.
\endproclaim

\demo{Proof}
We need to prove that for any element $e'\in\O$ there exists
an $\ssl_2$-triple $\langle e',h',f'\rangle$ such that
$h'\in i\k$, where $\k$ is a {\it fixed} compact real form of $\g$.
Clearly it is sufficient to prove that for an {\it arbitrary}
compact real form $\k$ there exists an $\ssl_2$-triple
$\langle e,h,f\rangle$ with $h\in i\k$.
According to the proof of the previous Lemma
we should choose 
$h$ to be a unique characteristic such that
$B(\overline h,\z^u_{\g}(e))=0$, where
$x\to\overline x$ denotes a complex conjugation in $\g$ with respect to 
the compact form $\k$.
It remains to prove that $h\in i\k$.
Since $B$ is a non-degenerate $\ad$-invariant scalar product on $\g$
and $\z_\g^u(e)=\n_+$
it follows that $\overline h\in\p$, where $\p=\g_0\oplus\n_+$.
Let $\l$ be some ``standard'' compact real form of $\g$
such that $h\in i\l$ and $\tilde\n_{\pm}=\n_{\mp}$,
where 
$x\to\tilde x$ denotes a complex conjugation in~$\g$ with respect to 
the compact form $\l$.
Let $P\subset G$ be a parabolic subgroup of $G$
with the Lie algebra $\p$, let $H\subset P$ be its Levi subgroup
with the Lie algebra $\g_0$.
Since all compact real forms of a semisimple complex Lie algebra
are conjugated by elements of any fixed Borel subgroup
it follows that there exists $g\in P$ such that
$\Ad(g)\k=\l$.
Therefore 
$$\widetilde{\Ad(g)h}=\Ad(g)\overline h\subset\Ad(g)(\p)=\p.$$
We can express $g$ as a product $uz$, where $u\in\exp(\n_+)$, $\Ad(z)h=h$.
Then 
$\widetilde{\Ad(g)h}=\widetilde{\Ad(u)h}$.
If $u$ is not the identity element of $G$
then $\Ad(u)h=h+\xi$, where $\xi\in\n_+$ and $\xi\ne0$.
Therefore $\widetilde{\Ad(u)h}=-h+\tilde\xi$.
But $\tilde\xi\in\n_{-}$ and hence
$\widetilde{\Ad(u)h}\not\in\p$, contradiction.
Therefore $u$ is trivial and since
$\Ad(z)h=h$ we finally get
$$\overline h=\tilde h=-h.\qed$$
\enddemo

Now we shall try to reverse this argument.

\proclaim{Lemma 1.3}
Suppose that $\O=\Ad(G)e$ is a Moore--Penrose orbit. Then $\z_\g^u(e)=\n_+$.
\endproclaim

\demo{Proof}
We choose a standard compact real form $\l$ 
as in the proof of the previous Lemma.
Clearly, $\z_\g^u(e)$ is a graded subalgebra of $\n_+=\mathop{\oplus}\limits_{k>0}\g_k$.
Suppose, on the contrary, that $\z_\g^u(e)\ne\n_+$. Let $\xi\in\g_p$, $p>0$,
be a homogeneous element that does not belong to $\z_\g^u(e)$.
Let $u=\exp(\xi)$. Let $e'=\Ad(u)e$. We claim that all characteristics
of $e'$ do not belong to $i\l$. Indeed, all characteristics
of $e'$ have a form $\Ad(u)h+\Ad(u)x$, where $x\in\z_\g^u(e)$.
Suppose that for some $x$ we have $\Ad(u)h+\Ad(u)x\in i\l$.
Since $h\in i\l$, $\tilde\n_{\pm}=\n_{\mp}$, and
$\Ad(u)(h+x)-h\in\n_+$ it follows that 
$\Ad(u)(h+x)=h$. In $\n_+$ modulo 
$\mathop{\oplus}\limits_{k>p}\g_k$ we obtain the equation
$[\xi,h]+x=0$, but $[h,\xi]=p\xi$ and therefore $\xi\in\z_\g^u(e)$.
Contradiction.
\qed\enddemo

\demo{Proof of Theorem~2}
Combining Lemma~1.2 and~1.3 we see that 
$\O$ is a Moore--Penrose orbit if and only if 
$\z_\g^u(e)=\n_+$.
It follows from the $\ssl_2$-theory that
$\dim\z_\g^u(e)=\dim\g_1+\dim\g_2$.
Therefore $\z_\g^u(e)=\n_+$ if and only if $\g_p=0$ for $p>2$.
Clearly, this is precisely equivalent to $\hht(\O)=2$.
In this case for any $e\in\O$ there exists a unique $\ssl_2$-triple
$\langle e,h,f\rangle$ such that $h\in i\k$ by Lemma~1.1.
\qed\enddemo

\head \S2. Moore--Penrose orbits in simple graded Lie algebras\endhead
In this section we shall prove Theorems~1 and~3. 
We shall use essentially the same arguments as in the proof
of Theorem~2, we shall only need to adapt them to the homogeneous situation.
First we shall prove
Theorem~3 and then deduce Theorem~1 from it.

Recall that $\g$ is a $\Z$-graded simple Lie algebra, 
$\g=\mathop{\oplus}\limits_{k\in\Z}\g_k$, 
$G$ is a corresponding simple simply-connected group.
$P\subset G$ is a parabolic subgroup with Lie algebra 
$\g=\mathop{\oplus}\limits_{k\ge0}\g_k$,
$L\subset P$ is a Levi subgroup with Lie algebra $\g_0$,
$\k_0$ is a compact real form of $\g_0$.

Suppose $\O=\Ad(L)e\subset\g_k$. Take any homogeneous $\ssl_2$-triple
$\langle e,h,f\rangle$. Then $h$ defines a grading 
$\g_0=\mathop{\oplus}\limits_{n\in\Z}\g_0^n$,
such that $\ad(h)|_{\g_0^n}=n\cdot\Id$. Denote $\mathop{\oplus}\limits_{n>0}\g_0^n$ by $\n_+$
and $\mathop{\oplus}\limits_{n<0}\g_0^n$ by $\n_-$.

\proclaim{Lemma 2.1}
If $\ad(e)\n_+=0$ then $\O=\Ad(L)e$ is a Moore--Penrose orbit.
\endproclaim

\demo{Proof}
We need to prove that for any element $e'\in\O$ there exists
a homogeneous $\ssl_2$-triple $\langle e',h',f'\rangle$ such that
$h'\in i\k_0$, where $\k_0$ is a {\it fixed} compact real form of~$\g_0$.
Clearly it is sufficient to prove that for an {\it arbitrary}
compact real form $\k_0$ of $\g_0$ there exists a homogeneous $\ssl_2$-triple
$\langle e,h,f\rangle$ with $h\in i\k_0$.
Let us start with an arbitrary homogeneous $\ssl_2$-triple
$\langle e,h,f\rangle$.
It is well-known that the space $\cal H$ of all possible characteristics
is an affine space $h+\z_\g^u(e)$.
Therefore the space of all homogeneous characteristics
is an affine space $h+\z_{\g_0}^u(e)=h+\n_+$.
Arguing as in the proof of Lemma~1.2,
let us change a characteristic $h$ in such a way that
$B(\overline h,\n_+)=0$, where
$x\to\overline x$ denotes a complex conjugation in $\g_0$ with respect to 
the compact form $\k_0$, and $B$ is Killing form in $\g$, not in $\g_0$!
It remains to prove that $h\in i\k_0$.
Since $B$ is a non-degenerate $\ad$-invariant scalar product on $\g_0$
it follows that $\overline h\in\q$, where $\q=\g_0^0\oplus\n_+$.
Let $\l_0$ be some ``standard'' compact real form of $\g_0$
such that $h\in i\l_0$ and $\tilde\n_{\pm}=\n_{\mp}$,
where 
$x\to\tilde x$ denotes a complex conjugation in $\g_0$ with respect to 
the compact form $\l_0$.
Let $Q\subset L$ be a parabolic subgroup of $L$
with Lie algebra $\q$, let $H\subset Q$ be its Levi subgroup
with Lie algebra $\g_0^0$.
There exists $g\in Q$ such that
$\Ad(g)\k_0=\l_0$.
(The conjugation theorem is usually stated only for semi-simple
Lie algebras, while $\g_0$ is only reductive. But according
to our conventions (see Remark in the Introduction),
the conjugation theorem holds for $\g_0$ as well.)
Therefore 
$$\widetilde{\Ad(g)h}=\Ad(g)\overline h\subset\Ad(g)(\q)=\q.$$
We can express $g$ as a product $uz$, where $u\in\exp(\n_+)$, $\Ad(z)h=h$.
Then 
$\widetilde{\Ad(g)h}=\widetilde{\Ad(u)h}$.
If $u$ is not the identity element of $G$
then $\Ad(u)h=h+\xi$, where $\xi\in\n_+$ and $\xi\ne0$.
Therefore $\widetilde{\Ad(u)h}=-h+\tilde\xi$.
But $\tilde\xi\in\n_{-}$ and hence
$\widetilde{\Ad(u)h}\not\in\q$, contradiction.
Therefore $u$ is trivial and since
$\Ad(z)h=h$ we finally get
$$\overline h=\tilde h=-h.\qed$$
\enddemo

\proclaim{Lemma 2.2}
Suppose that $\O=\Ad(L)e$ is a Moore--Penrose orbit. 
Then $\ad(e)\n_+=0$.
\endproclaim

\demo{Proof}
The proof is parallel to the proof of Lemma~1.3.
%We choose a standard compact real form $\l_0$ 
%as in the proof of the previous Lemma.
%Clearly, $\z_{\g_0}^u(e)$ is a graded subalgebra of $\n_+=\mathop{\oplus}\limits_{k>0}\g_0^k$.
%Suppose, on the contrary, that $\z_{\g_0}^u(e)\ne\n_+$. 
%Let $\xi\in\g_0^p$, $p>0$,
%be a homogeneous element that does not belong to $\z_{\g_0}^u(e)$.
%Let $u=\exp(\xi)$. Let $e'=\Ad(u)e$. We claim that all characteristics
%of $e'$ don't belong to $i\l_0$. Indeed, all characteristics
%of $e'$ have a form $\Ad(u)h+\Ad(u)x$, where $x\in\z_{\g_0}^u(e)$.
%Suppose that for some $x$ we have $\Ad(u)h+\Ad(u)x\in i\l_0$.
%Since $h\in i\l_0$, $\tilde\n_{\pm}=\n_{\mp}$, and
%$\Ad(u)(h+x)-h\in\n_+$ it follows that 
%$\Ad(u)(h+x)=h$. In $\n_+$ modulo 
%$\mathop{\oplus}\limits_{k>p}\g_0^k$ we obtain the equation
%$[\xi,h]+x=0$, but $[h,\xi]=p\xi$ and therefore $\xi\in\z_{\g_0}^u(e)$.
%Cotradiction.
\qed\enddemo

\demo{Proof of Theorem~3}
Combining Lemma~2.1 and~2.2 we see that 
$\O$ is a Moore--Penrose orbit if and only if 
$\ad(e)\g_0^n=0$ for any $n>0$.
In this case for any $e'\in\O$ there exists a unique homogeneous 
$\ssl_2$-triple $\langle e',h',f'\rangle$ such that
$h'\in i\k_0$ by Lemma~1.1.
\qed\enddemo

\demo{Proof of Theorem~1}
We should show that the condition
$\ad(e)\g_0^n=0$ for any $n>0$ is satisfied always
if $\g_k=0$ for $|k|>1$.
Suppose that $x\in\g_0^n$, $n>0$.
If $\ad(e)x\ne0$ then there exists an element $y\in\g_1$
such that $\ad(h)y=(n+2)y$.
Since $\ad(e)\g_i\subset\g_{i+1}$ it follows from $\ssl_2$-theory
that there exists a non-zero element $z\in\g_{1-(n+2)}$.
But $1-(n+2)<-1$. Contradiction.
\qed\enddemo

\head \S3. Jordan pairs\endhead

A Jordan pair is a pair of vector spaces $(V_+, V_-)$
with trilinear multiplications
$$V_\pm\otimes V_\mp\otimes V_\pm\to V_\pm,\quad x\otimes y\otimes z\mapsto \{xyz\},$$
which satisfy a certain set of axioms (see \cite{Lo}).
Fortunately, there is no need to write them down due to the
following fundamental observation.
$V_+$ and $V_-$ form a Jordan a triple if and only if
there exists a short graded Lie algebra $\g_{-1}\oplus\g_0\oplus\g_1$
such that 
$$\g_{-1}=V_-,\quad \g_1=V_+,\quad\text{and}\ \{x,y,z\}={1\over2}[[x,y],z].$$
In fact, this construction provides a bijection of 
the set of Jordan pairs up to an isomorphism and the set
of short graded Lie algebras up to a certain equivalence relation.
We refer the reader to \cite{Lo} and \cite{Ja} for more details
about Jordan pairs and Jordan algebras used throughout this section.

\example{Example 1}
Take $V_+=\Mat_{n,m}$, $V_-=\Mat_{m,n}$. Then $(V_+, V_-)$ is a Jordan pair
with respect to trilinear maps $\{ABC\}={1\over2}(ABC+CBA)$ (matrix multiplication).
This Jordan pair corresponds to a short grading of $\ssl_{n+m}$.
\endexample

\example{Example 2}
Suppose that $A$ is a Jordan algebra, that is, 
an algebra with a unit such that the bilinear multiplication 
in $A$ satisfies two axioms
$$ab=ba\ \text{(commutativity)},\quad ((aa)b)a=(aa)(ba)\ \text{(Jordan axiom)}.$$
For example, we can take any associative algebra and define a new multiplication
by a formula $a*b={1\over2}(ab+ba)$. 
This will be a (special) Jordan algebra.
Any Jordan algebra $A$ corresponds to a Jordan pair $(V_+, V_-)$
defined as follows:
$$V_+=V_-=A, \quad \{abc\}=(ab)c+(bc)a-(ac)b\quad \text{(Jordan triple product)}.$$
Simple Jordan algebras correspond to shortly graded simple Lie algebras
such that the corresponding Hermite homogeneous space $G/P$
(recall that $P$ is a parabolic subgroup with a Lie algebra $\g_0\oplus\g_1$)
has a tube type, or, equivalently, if there exists
an $L$-invariant hypersurface in~$\g_1$, where $L$ is a Levi subgroup of $P$.
\endexample

We shall be interested only in Jordan pairs arising from
shortly graded complex simple Lie algebras. It can be shown
that these Jordan pairs are precisely simple complex Jordan pairs.
For simplicity we shall use the term `Jordan pairs' 
only for these pairs.

In the Introduction we defined a MP-inverse for shortly graded
simple Lie algebras. In the language of Jordan pairs
a MP-inverse is some map $V_\pm\to V_\mp$.
The first aim of this section is to define a MP-inverse
entirely in terms of trilinear maps $\{\cdot,\cdot,\cdot\}$.
First let us give some definitions and lemmas.
For any Jordan pair $(V_+,V_-)$ we denote the corresponding shortly graded
simple Lie algebra by~$\g=\g_{-1}\oplus\g_0\oplus\g_1$.

\definition{Definition}
{\it A Killing pairing} $B(\cdot,\cdot)$ of a Jordan pair
is a bilinear map 
$$V_\pm\otimes V_\mp\to\C\ \text{given by}\ 
x\otimes y\mapsto\Tr\{x,y,\cdot\}.$$
\enddefinition

\example{Example}
Suppose that the Jordan pair corresponds to a Jordan algebra A.
Then 
$$B(x,y)=\Tr\{x,y,\cdot\}=
\Tr\big\{(ab)\cdot+(b\cdot)a-(a\cdot)b\big\}=
\Tr\big\{(ab)\cdot\big\}+\Tr\big\{a(b\cdot)-b(a\cdot)\big\}=
\Tr\big\{(ab)\cdot\big\}.$$
This is a usual definition of a scalar product in a Jordan algebra
(up to a positive multiple).
\endexample

\proclaim{Proposition 3.1}
A Killing pairing is symmetric and non-degenerate.
It coincides up to a positive multiple 
with a restriction of a Killing form of $\g$
on $\g_{-1}\oplus\g_1$.
\endproclaim

\demo{Proof}
Since $\g_{-1}$ and $\g_1$ are dual $\g_0$-modules it follows that
for any $\xi\in\g_0$ we have 
$$\Tr [\xi,\cdot]|_{\g_{-1}}=-\Tr [\xi,\cdot]|_{\g_{1}}.$$
Therefore,
for any $x\in V_+$, $y\in V_-$, we get
$$B(x,y)=\Tr\{x,y,\cdot\}={1\over2}\Tr[[x,y],\cdot]|_{\g_{1}}=
-{1\over2}\Tr[[x,y],\cdot]|_{\g_{-1}}=
{1\over2}\Tr[[y,x],\cdot]|_{\g_{-1}}=B(y,x).$$
This proves symmetry. To show that $B$ is non-degenerate it is sufficient
to prove that
$B(\cdot,\cdot)$ coincides up to a positive multiple 
with a restriction of a Killing form $(\cdot,\cdot)$ of $\g$
on $\g_{-1}\oplus\g_1$.
We choose a Cartan subalgebra $\h\subset\g_0$.
Let $\Delta$ be the corresponding root system.
For any root $\alpha\in\Delta$ we choose a root vector $e_\alpha\in\g$
in such a way that $[e_\alpha,e_{-\alpha}]=h_\alpha$, where $h_\alpha\in\h$
is a coroot of $\alpha$, so for any $\beta\in\Delta$ we have
$\beta(h_\alpha)={2(\beta,\alpha)\over (\alpha,\alpha)}$.
We have a decomposition $\Delta=\Delta_{-1}\cup\Delta_0\cup\Delta_1$,
where $\alpha\in\Delta_k$ if and only if $e_\alpha\in\g_k$.
We also take a set of simple roots $\Pi_0\subset\Delta_0$
and a root $\gamma\in\Delta_1$ such that the system $\Pi_0\cup\{\gamma\}$
is a set of simple roots for $\Delta$.
Let $\Delta_0^+$ be a set of positive roots of $\Delta_0$ corresponding to $\Pi_0$.
Then $\Delta^+=\Delta_0^+\cap\Delta_1$ is a set of positive roots for $\Delta$.
It suffices to show that for any $\alpha\in\g_1$, $\beta\in\g_{-1}$
we have $B(e_\alpha,e_\beta)=c(e_\alpha,e_\beta)$, where $c>0$
does not depend on $\alpha$ and $\beta$.
If $\alpha+\beta\ne0$ then clearly $(e_\alpha,e_\beta)=0$.
But in this $B(e_\alpha,e_\beta)$ is also equal to zero,
because $[[e_\alpha,e\beta],\cdot]$ is a nilpotent operator.
Suppose now that $\beta=-\alpha$.
Then $(e_\alpha,e_{-\alpha})={2\over(\alpha,\alpha)}$.
On the other hand,
$$B(e_\alpha,e_{-\alpha})=
{1\over2}\Tr[h_\alpha,\cdot]|_{\g_{1}}=
\sum_{\beta\in\Delta_1}{2(\beta,\alpha)\over (\alpha,\alpha)}=
{2(\rho_1,\alpha)\over (\alpha,\alpha)},\ \text{where}
\ \rho_1=\sum\limits_{\beta\in\Delta_1}\beta.$$
So it suffices to prove that $(\rho_1,\alpha)$ is positive
and does not depend on the choice of $\alpha\in\Delta_1$.
In fact, the first claim will follow from the second, because
then $(\rho_1,\alpha)={1\over\#\Delta_1}(\rho_1,\rho_1)>0$.
So let us prove the second claim.
Any root $\alpha\in\Delta_1$ is a positive linear combination of $\gamma$
and some simple roots from $\Pi_0$.
Therefore we need to prove that for any $\delta\in\Pi_0$ we have
$(\rho_1,\delta)=0$. But 
$$\rho_1=\rho-\rho_0,\ \text{where}\ 
\rho=\sum\limits_{\beta\in\Delta^+}\beta,
\ \rho_0=\sum\limits_{\beta\in\Delta_0^+}\beta.$$
Therefore 
$$(\rho_1,\delta)=(\rho,\delta)-(\rho_0,\delta)=
(\delta,\delta)-(\delta,\delta)=0.\qed$$
\enddemo

\definition{Definition}
A pair $\omega$ of antilinear maps $V_\pm\to V_\mp$ 
is called a Cartan involution of a Jordan pair if 
$$\omega^2=\Id,\quad \{\omega(x),\omega(y),\omega(z)\}=\omega\{x,y,z\},$$
$$\text{Hermite form}\ H(x)=B(x,\omega(x))\ \text{is positive definite on}
\ V_\pm.$$
\enddefinition

Suppose that $\tilde\omega$ is a Cartan involution of $\g$ 
(so $\g^{\tilde\omega}$ is a compact real form of $\g$)
such that $\tilde\omega(\g_k)=\g_{-k}$.
Let $\sigma\in\Aut(\g)$ be defined as follows: 
$$\sigma|_{\g_0}=\Id,\quad\sigma|_{\g_{-1}\oplus\g_1}=-\Id.$$
Clearly, $\hat\omega=\sigma(\tilde\omega)$ is an antilinear involution of $\g$.
Then $\hat\omega$
restricted to $\g_{-1}\oplus\g_1$ is a Cartan involution $\omega$
of a corresponding Jordan pair.
In particular, any Jordan pair has a Cartan involution.

\proclaim{Proposition 3.2}
The correspondence $\tilde\omega\to\omega$ is bijective.
\endproclaim

\demo{Proof}
The set of commutators $[x,y]$ for $x\in\g_1$, $y\in\g_{-1}$
spans $\g_0$. Since $\tilde\omega([x,y])=[\omega(x),\omega(y)]$
it follows that this correspondence is injective.
Suppose now that $\omega$ is a Cartan involution of a Jordan pair.
Then we define $\tilde\omega$ on $\g_{-1}\oplus\g_1$ as $-\omega$ and
we define $\tilde\omega$ on $\g_0$ by setting
$\tilde\omega([x,y])=[\omega(x),\omega(y)]$.
Let us show that $\tilde\omega$ is well-defined.
Suppose that $\xi=[x_1,y_1]+\ldots+[x_k,y_k]=0$
for $x_i\in\g_1$, $y_i\in\g_{-1}$.
We need to show that
$\xi_\omega=[\omega(x_1),\omega(y_1)]+\ldots+[\omega(x_k),\omega(y_k)]=0$.
The representation of $\g_0$ on $\g_1$ is faithful.
Therefore it suffices to prove that for any $z\in\g_{-1}$ we have
$[\xi_\omega,\omega(z)]=0$.
But
$$[\xi_\omega,\omega(z)]=2\sum_{i=1}^k\{\omega(x_i),\omega(y_i),\omega(z)\}=
2\sum_{i=1}^k\omega\{x_i,y_i,z\}=\omega[\xi,z]=0.$$
Therefore $\tilde\omega$ is well-defined.
Clearly $\tilde\omega^2=\Id$ and 
$[\tilde\omega(x),\tilde\omega(y)]=\tilde\omega[x,y]$.
It remains to prove that the Hermite form 
$\tilde H(x)=-(x,\tilde\omega(x))$
is positive definite on $\g$. Since the decomposition
$\g=\g_{-1}\oplus\g_0\oplus\g_1$ is orthogonal with respect to $\tilde H$
and the restriction of $\tilde H$ on $\g_{-1}$ and $\g_1$
coincides with $H$ up to a positive multiple,
$\tilde H$ is positive definite on $\g_{-1}\oplus\g_1$.
There exists a Cartan involution $\tau$ compatible with $\tilde\omega$
and such that $\tau(\g_0)=\g_0$ (see \cite{VO}).
Then, clearly, $\tau(\g_{\pm1})=\g_{\mp1}$.
Therefore, $\tau\tilde\omega$ is an involution preserving the grading.
We need to show that $\tau\tilde\omega=\Id$.
It is sufficient to show that $\tau\tilde\omega|_{\g_1}=\Id$.
Suppose, on the contrary, that there exists $x\in \g_1$, $x\ne0$, such that
$\tau\tilde\omega(x)=-x$. Then
$$0>(x,\tilde\omega(x))=(\tau^2(x),\tilde\omega(x))=
(\tau(x), \tau\tilde\omega(x))=-(x, \tau(x))>0.$$
Contradiction.
\qed\enddemo

Now we are ready to give a new definition of a Moore--Penrose inverse.

\definition{Definition}
Suppose that $(V_+,V_-)$ is a Jordan pair.
Then for any $A\in V_\pm$ its MP-inverse is an element $A^+\in V_\mp$
such that
$$\{AA^+A\}=A,\quad \{A^+AA^+\}=A^+, \eqno(*)$$
$$\{AA^+\cdot\},\ \{A^+A\cdot\}
\ \text{are Hermite operators with respect to}\ H.\eqno(**)$$
\enddefinition

In fact we shall see later that if one of operators in $(**)$ is Hermite
then another one is a Hermite operator automatically.

\example{Example}
Suppose that the Jordan pair corresponds to a Jordan algebra A.
An element $b\in A$ is called a (usual) inverse of an element $a$
if $ab=1$, $(aa)b=a$. This definition does not look very symmetric,
but in fact it is. So automatically $(bb)a=b$.
Let us show that in this case $b$
coincides with a MP-inverse $a^+$.
Indeed, $\{aba\}=(ab)a+(ba)a-(aa)b=a+a-a=a$.
Similarly $\{bab\}=b$. Moreover, one can show that operators $(a\cdot)$
and $(b\cdot)$ commute (see \cite{Ja}), and, therefore,
$\{ab\cdot\}$ and $\{ba\cdot\}$ are identity operators and, hence, 
Hermite operators with respect to $H$.
\endexample

\proclaim{Proposition 3.3}
New definition of a MP-inverse coincides with given in the Introduction.
In particular, MP-inverse in Jordan pairs exists and is unique.
\endproclaim

\demo{Proof}
Indeed, conditions $(*)$ mean that $\langle A, [A,A^+], A^+\rangle$
is a homogeneous $\ssl_2$-triple in $\g$. We choose a compact form $\k_0$
in $\g_0$ such that $\k_0=\g_0^{\tilde\omega}$.
Since the representations of $\g_0$ in $\g_1$ and $\g_{-1}$ are faithful,
any of operators in $(**)$ is Hermite if and only if $[A,A^+]\in i\k_0$.
\qed\enddemo

Now we can describe a MP-inverse arising from short gradings of 
$E_6$ and $E_7$.
We start with~$E_7$.
Let $\Ca$ denote the algebra of split Cayley numbers over $\C$.
Let $\OO\subset\Ca$ be the Cayley division algebra of octonions.
Let $\A$ be the Albert algebra of Hermite $(3\times3)$-matrices
over $\Ca$ (with respect to the canonical anti-involution
in $\Ca$). This is a complex Jordan algebra with respect to the
symmetrisation of matrix multiplication.
The corresponding Jordan pair arises from the short grading of $E_7$.
Let $\A(\OO)\subset\A$ be the real subalgebra of matrices with entries in~$\OO$.
It is well-known that the scalar product $\Tr((ab)\cdot)$
is positive definite on $\A(\OO)$. 
Therefore the complex conjugation $\omega$ of $\A$ with respect to $\A(\OO)$
is a Cartan involution and we have all the information necessary to write
down equations $(*,**)$.

The Jordan pair corresponding to the short grading of $E_6$
is, in fact, `a subpair' of a previous one.
Namely, $V_+=V_-=\Ca\oplus\Ca$. We can consider both of them as matrices
from $\A$ of the form
$$\left(\matrix 0&c_1&c_2\cr \overline c_1&0&0\cr \overline c_2&0&0\cr
\endmatrix\right).$$
Then the Jordan triple product in $\A$ defines the trilinear
maps for this Jordan pair.
So the Moore--Penrose inverse in this case is the restriction
of a previous MP-inverse.

\head \S4. Moore--Penrose parabolic subgroups\endhead

The definition of Moore--Penrose parabolic 
subgroups given in the Introduction does not
look very natural because it depends on the choice of a grading.
Let us give a more transparent equivalent definition.
Suppose that $G$ is a simple simply--connected Lie group,
$P\subset G$ is a parabolic subgroup, $L\subset P$ is a Levi subgroup,
$Z\subset L$ is a connected component of its center.
Then a Lie algebra $\g$ of $G$ admits a natural $\Z^m$-grading,
where $\Z^m$ is a character group of $Z$.
Clearly a Lie algebra $\l$ of $L$ is just a zero component
of this grading. 
In fact any $\Z$-grading of $\g$ ``compatible'' with $P$
can be obtained from this $\Z^m$-grading via some homomorphism $\Z^m\to \Z$.
Moreover, for generic homomorphism non-zero homogeneous components
will be the same for $\Z$- and $\Z^m$-grading.
Let $\k_0$ be a compact real form of $\l$.
It is easy to see that the definition of Moore--Penrose parabolic
subgroups given in the Introduction is equivalent
to the following:

\definition{Definition}
$P$ is called {\it a Moore--Penrose parabolic subgroup\/}
if for any $e\in\g_\alpha$, $\alpha\in\Z^m$, $\alpha\ne0$,
there exists a $\Z^m$-homogeneous $\ssl_2$-triple $\langle e,h,f\rangle$
such that $h\in i\k_0$. In this case $f$ is called a MP-inverse of $e$.
\enddefinition

For any $\alpha\in\Z^m$, $\alpha\ne0$,
let us consider a reductive subalgebra
$\g^\alpha=\mathop{\oplus}\limits_{\beta\in\R\alpha}\g_\beta$.
Then $\g^\alpha$ is a Levi subalgebra of $\g$
and, clearly, $P$ is a Moore--Penrose parabolic subgroup
if and only if the maximal parabolic subalgebra 
$\mathop{\oplus}\limits_{\beta\in\R_{\ge0}\alpha}\g_\beta$ of $\g^\alpha$
is a Moore--Penrose parabolic subalgebra for any~$\alpha$.
In particular, in order to describe all Moore--Penrose 
parabolic subgroups of $G$ it suffices to find
all Moore--Penrose maximal parabolic subgroups
in simple components of Levi subgroups of $G$.

\demo{Proof of Theorem~4}
All simple components of Levi subgroups of $\SL_n$
are again simple groups of type $\SL_k$. 
All maximal parabolic subgroups in $\SL_k$
have aura and therefore are Moore--Penrose parabolic subgroups
by Theorem~1.
\qed\enddemo

In the rest part of this section we shall find all
Moore--Penrose maximal parabolic subgroups in remaining
classical groups $\SO_n$ and $\Sp_n$.
Let $G$ be one of this groups, $P$ be its maximal parabolic subgroup.
We denote their Lie algebras by $\g$ and $\p$.
The corresponding grading has a form
$\g=\mathop{\oplus}\limits_{k\in\Z}\g_k$. 
Recall that $L$ is a Levi subgroup of $P$ such that $\g_0$ is
a Lie algebra of $L$.
It easy to see that
$\g_k=0$ either for $|k|>1$ or for $|k|>2$.
In the first case $P$ has an aura and therefore 
is Moore-Penrose by Theorem~1. 
We take a usual bijection between maximal parabolic subgroups
and simple roots of the corresponding algebra.
In Bourbaki-numbering
of simple roots $P$ has an aura if and only if $P$ corresponds
to one of simple roots $\alpha_1$, $\alpha_n$ ($B_n$-case);
$\alpha_n$ ($C_n$-case); $\alpha_1$, $\alpha_{n-1}$, $\alpha_n$ ($D_n$-case).
Now let us consider other possibilities.
Clearly all $L$-orbits in~$\g_2$ are Moore--Penrose 
(see Example~1 in the Introduction).
Now we shall classify all Moore--Penrose $L$-orbits in $\g_1$.
But first we shall reformulate the problem in linear--algebraic terms.

Suppose that $U=\C^k$, $V=\C^n$ are complex vector spaces
with standard Hermite scalar products. In the orthogonal case
we suppose that $k\ge2$, $n\ge3$. In the symplectic case
we assume that $k\ge1$, $n\ge2$, and $n$ is even.
We choose a symmetric (resp.~skew-symmetric) $2$-form $\omega$ in $V$
with matrix $I=E$ (resp.~$I=\left(\matrix 0&E\cr-E&0\endmatrix\right)$),
where $E$ is an identity matrix. Let 
$\Ow$ be a special orthogonal (resp.~symplectic) group
corresponding to $\omega$. Let $\ow$ be a corresponding Lie algebra.
Then $L$, $\g_0$, $\g_1$, $\g_{-1}$ have
the following interpretation:
$$L=\SL(U)\times\Ow\times\C^*,\quad
\g_0=\ssl(U)\oplus\ow\oplus\C,$$
$$\g_1=\Hom(U,V),\quad \g_{-1}=\Hom(V,U).$$
The action of $\SL(U)\times\Ow$ on $\Hom(U,V)$ is standard,
$\C^*$ acts by homotheties. So the action of $L$ on $\g_1$
has finitely many orbits $\O(a,b)$ indexed by $a=\rk(F)$ and 
$b=\dim\Ker\,\omega|_{\Im F}$
for $F\in\Hom(U,V)$. For any $F\in\Hom(U,V)$ its Moore--Penrose
inverse (if exists) is defined as a unique $G\in\Hom(V,U)$
such that
$$GF,\  FG-(FG)^\#\ \text{are Hermite operators},\eqno(*)$$
$$F=2FGF-(FG)^\#F,\quad G=2GFG-G(FG)^\#,\eqno(**)$$
where for any $A\in\Hom(V,V)$ we denote by $A^\#$
its adjoint operator with respect to $\omega$.
Now let us prove the following proposition

\proclaim{Proposition}
An $L$-orbit $\O(a,b)\subset\g_1$ is Moore--Penrose
if and only either $b=0$ or $a=b$.
\endproclaim

\demo{Proof}
Suppose first that $b=0$, $F\in\O(a,b)$.
Then the restriction of $\omega$ on $\Im(F)$ is non degenerate.
Let $\Im(F)^\perp$ denote its orthogonal complement with respect to $\omega$.
Let $\Ker(F)^\perp$ denote the orthogonal complement of $\Ker(F)$
with respect to the Hermite form on $U$.
Let $\tilde F\in\Hom(\Ker(F)^\perp, \Im(F)$ be an operator induced by $F$.
Let $G\in\Hom(V,U)$ be an operator defined as follows:
$G|_{\Im(F)^\perp}=0$, $G|_{\Im(F)}=\tilde F^{-1}$.
Then $GF$ is a Hermite projector on $\Ker(F)^\perp$.
Since $FG$ is an orthogonal projector on $\Im(F)$ with respect to $\omega$,
$(FG)^\#=FG$, therefore $FG-(FG)^\#=0$ is a Hermite operator,
now relations $(**)$  are obvious.

Suppose now that $b=a$, $F\in\O(a,b)$.
Then the restriction of $\omega$ on $\Im(F)$ is equal to~$0$.
Let $\Im(F)^0=I\overline{\Im(F)}$ (recall that $I$ is a matrix of $\omega$,
bar denotes the complex conjugation).
Then $\Im(F)\cap\Im(F)^0=0$, the restriction of $\omega$
on $\Im(F)\oplus\Im(F)^0$ is non--degenerate and the orthogonal
complement $V'$ of $\Im(F)\oplus\Im(F)^0$
with respect to $\omega$ coincides with the orthogonal complement
of $\Im(F)\oplus\Im(F)^0$ with respect to a Hermite form.
Let $\Ker(F)^\perp$ denote the orthogonal complement of $\Ker(F)$
with respect to the Hermite form on $U$.
Let $\tilde F\in\Hom(\Ker(F)^\perp, \Im(F)$ be an operator induced by $F$.
Let $G\in\Hom(V,U)$ be an operator defined as follows:
$G|_{V'}=G|_{\Im(F)^0}=0$, $G|_{\Im(F)}={1\over2}\tilde F^{-1}$.
Then $GF$ is a Hermite projector on $\Ker(F)^\perp$.
It is clear that $FG$ is equal to~$0$ on $V'$ and $\Im(F)^\circ$
and is an identity operator on $\Im(F)$.
Therefore its adjoint operator $(FG)^\#$
is equal to~$0$ on $V'$ and $\Im(F)$
and is a minus identity operator on $\Im(F)^0$.
Therefore $FG-(FG)^\#$ is an orthogonal (and Hermite)
projector on $\Im(F)\oplus\Im(F)^0$.
Now relations $(**)$ are obvious.

It remains to prove that if $0<b<a$ then $\O(a,b)$
is not a Moore--Penrose orbit.
Choose a subspace $L\subset V$ such that $\dim\Ker\,\omega|_L=b$,
let $L_0=\Ker\,\omega|_L$. 
Let $U=U_0\oplus U_1\oplus U_2$
be an orthogonal (with respect to the Hermite form)
direct sum of subspaces such that $\dim U_0=b$, $\dim U_0+\dim U_1=a$.
Let $F\in\Hom(U,V)$ be a linear operator such that
$F|_{U_2}=0$, $F(U_0)=L_0$, $F(U_0\oplus U_1)=L$,
$F(U_1)$ is not orthogonal to $L_0$ with respect to the Hermite form.
We set $L_1=F(U_1)$.
We claim that $F$ does not have a Moore--Penrose inverse.
Suppose, on the contrary, that $G$ is a Moore--Penrose inverse of $F$.
Since $GF$ is Hermite, we see that $G(L)\subset U_0\oplus U_1$.
If $v\in L_0$, $v'\in V$ then
$\omega((FG)^\#v,v')=\omega(v,FGv')=0$, because $\omega(L_0,L)=0$.
Therefore, $(FG)^\#|_{L_0}=0$. It follows that $F|_{U_0}=2FGF|_{U_0}$.
So, $FG|_{L_0}={1\over2}E$, $GF|_{U_0}={1\over2}E$.
Since $GF$ is Hermite, it should preserve $U_1$, therefore
$G|_{L_1}\subset U_1$. Finally, we see that $FG$ preserves both
$L_0$ and $L_1$. To obtain a contradiction it suffices to show
that $FG$ is a Hermite operator on $L$ (because we know that $L_0$
and $L_1$ are not orthogonal).
But since $F=2FGF-(FG)^\#F$, it follows that
$2FG-(FG)^\#=E$ on $L$. Therefore $FG=E-(FG-(FG)^\#)$ is Hermite,
because we know that $FG-(FG)^\#$ is Hermite.
\qed\enddemo

It is clear that all $L$-orbits are Moore--Penrose
(under our restrictions on $k$ and $n$)
if and only if $\omega$ is symplectic and either $k\le2$
or $n=2$. Therefore we have a following Corollary

\proclaim{Corollary}
1) All Moore--Penrose maximal parabolic subgroups in $\SO(n)$
have aura.\par
2) Moore--Penrose maximal parabolic subgroups without aura 
in $\Sp(2n)$ correspond to simple roots $\alpha_1$, $\alpha_2$, $\alpha_{n-1}$.
\endproclaim

\head \S5. Moore--Penrose inversion and rank stratification\endhead
In this section we prove Theorem~5.
Let us start with some lemmas.

Suppose that $G$ is a connected reductive group with a Lie algebra $\g$.
For any elements $x_1,\ldots,x_r\in\g$ let 
$\langle x_1,\ldots,x_r\rangle_{alg}$ denote
the minimal algebraic Lie subalgebra of $\g$ that contains
$x_1,\ldots,x_r$. By a theorem of Richardson \cite{Ri1}
$\langle x_1,\ldots,x_r\rangle_{alg}$ is reductive if and only
if an orbit of the $r$-tuple $(x_1,\ldots,x_r)$ in $\g^r$
is closed with respect to the diagonal action of~$G$.
Suppose now that $h_1,\ldots,h_r$ are semi-simple elements
of $\g$. Consider the closed variety 
$\hat\O=(\Ad(G)h_1,\ldots,\Ad(G)h_r)\subset\g^r$.
For any closed $G$-orbit $\O\subset\hat\O$ let
us denote by $G(\O)$ the conjugacy class of the reductive subalgebra
$\langle x_1,\ldots,x_r\rangle_{alg}$ for $(x_1,\ldots,x_r)\in\O$.

\proclaim{Lemma 5.1}
There are only finitely many conjugacy classes $G(\O)$.
\endproclaim

\demo{Proof}
We shall use induction on $\dim\g$. Suppose that the claim of Lemma~5.1
is true for all reductive groups $H$ with $\dim H<\dim G$.
Let $\z\subset\g$ be the center of $\g$, $\g'\subset\g$ be its
derived algebra. Consider two canonical homomorphisms
$$\g\mathop{\longrightarrow}\limits^\pi\g'\quad\text{and}\quad
\g\mathop{\longrightarrow}\limits^{\pi'}\z.$$
We take any closed $G$-orbit $\O\subset\hat\O$.
Let $(x_1,\ldots,x_r)\in\O$, $y_i=\pi(x_i)$ for $i=1,\ldots,r$.
Then 
$\langle y_1,\ldots,y_r\rangle_{alg}=
\pi(
\langle x_1,\ldots,x_r\rangle_{alg})$ and, therefore, is reductive.
Let us consider two cases.

Suppose first, that 
$\langle y_1,\ldots,y_r\rangle_{alg}=\g'$.
Then $\g'$ is a derived algebra of
$\langle x_1,\ldots,x_r\rangle_{alg}$
and, therefore, 
$\langle x_1,\ldots,x_r\rangle_{alg}=
\langle \pi'(h_1),\ldots,\pi'(h_r)\rangle\oplus\g'$.
In this case we get one conjugacy class.

Suppose now, that $\langle y_1,\ldots,y_r\rangle_{alg}\ne\g'$.
Then $\langle y_1,\ldots,y_r\rangle_{alg}$ is contained in some maximal
reductive Lie subalgebra of $\g'$. It is well-known (and not difficult
to prove) that in a semisimple Lie algebra there are
only finitely many conjugacy classes of maximal reductive
subalgebras. Let $\h'$ be one of them, $\h=\z\oplus\h'\subset\g$.
Let $H$ be a corresponding reductive subgroup of $G$.
It is sufficient to prove that for any closed $G$-orbit $\O$ of 
$\hat\O$ that meets $\h^r$ there are only finitely many 
possibilities for $G(\O)$. 
It easily follows from Richardson's Lemma \cite{Ri}
that for any $i$ the intersection $\Ad(G)h_i\cap\h$ is a union of finitely many closed
$H$-orbits, say $\Ad(H)h_i^1,\ldots,\Ad(H)h_i^{s_i}$.
It remains to prove that if for some $r$-tuple
$(x_1,\ldots,x_r)\in\Ad(H)h_1^{k_1}\times\ldots\Ad(H)h_r^{k_r}$
the corresponding subalgebra
$\langle x_1,\ldots,x_r\rangle_{alg}$ is reductive
then there are only finitely many possibilities 
for its conjugacy class. But this is precisely the
claim of Lemma for the group~H, which is true by the induction
hypothesis.
\qed\enddemo

Suppose that $\k$ is a compact real form of $\g$.

\proclaim{Lemma 5.1}
If $r$-tuple $(x_1,\ldots,x_r)$ belongs to $(i\k)^r$,
then its $G$-orbit is closed in $\g^r$.
\endproclaim

\demo{Proof}
Indeed, let $B$ be a non-degenerate $\ad$-invariant
scalar product on $\g$, which is negative-definite on $\k$.
Let $H(x)=-B(\overline x, x)$ be a positive-definite
$\k$-invariant Hermite quadratic form on $\g$,
where the complex conjugation is taken with respect to $\k$.
Let $H^r$ be a corresponding Hermite quadratic form on $\g^r$.
More precisely, $H^r(x_1,\ldots,x_r)=
H(x_1)+\ldots+H(x_r)$.
By a Kempf--Ness criterion \cite{PV}
in order to prove that the $G$-orbit 
of $(x_1,\ldots,x_r)$ is closed it is sufficient to prove 
that the real function $H^r(\cdot)$ has a critical point
on this orbit. Let us show that $(x_1,\ldots,x_r)$
is this critical point.
Indeed, for any $g\in\g$
$$-B(\overline x_1,[g,x_1])-\ldots-B(\overline x_r,[g,x_r])=
B(x_1,[g,x_1])+\ldots+B(x_r,[g,x_r])=0.\qed$$
\enddemo

Now let $G$ be a simple simply-connected Lie group,
let $\g$ be its Lie algebra with a $\Z$-grading
$\g=\mathop{\oplus}\limits_{k\in\Z}\g_k$.
Let $r$ be a maximal integer such that $\g_r\ne0$.
We are going to change slightly our habits and denote the non-positive part of the grading
$\mathop{\oplus}\limits_{k\le0}\g_k$ by $\p$. Let $P\subset G$
be a parabolic subgroup with the Lie algebra $\p$.
We shall identify $\g/\p$ with $\mathop{\oplus}\limits_{k>0}\g_k$.
Let $L\subset G$ be a connected reductive subgroup with Lie algebra~$\g_0$.
Let $V$ be an irreducible $G$-module.
If we choose a Cartan subalgebra $\t\subset\g_0$
then the grading of $\g$ originates from some $\Z$-grading
on $\t^*$. Therefore, there exists a 
$\Z$-grading $V=\mathop{\oplus}\limits_{k\in\Z}V_k$
such that $\g_iV_j\subset V_{i+j}$.
Let $R$ be a maximal integer such that $V_R\ne0$.
It is easy to see that $M_V=\mathop{\oplus}\limits_{k<R}V_k$
(notice that $M_R$ is an irreducible $L$-module).
Now we shall prove Theorem~5.

\demo{Proof of Theorem~5}
It is sufficient to prove that there exists a finite set of points
$\{x_1,\ldots,x_N\}\subset\g/\p$
such that for any $x\in\g/\p$ and for any $V$
we have $\rk\,\Psi_V(x)=\rk\,\Psi_V(x_i)$ for some $i$.
Recall that $L$ has finitely many orbits on each $\g_k$ \cite{Ri,Vi}.
We pick some $L$-orbit $\O_i$ in each $\g_i$.
Then it is sufficient to find a finite set of points as above
only for points $x\in\g/\p$ of a form $x=x_1+\ldots+x_r$, where $x_i\in\O_i$.
For any orbit $\O_i$ let ${\cal H}_i$ denote the set of all possible
homogeneous characteristics of all elements from $\O_i$.
Clearly ${\cal H}_i$ is a closed $\Ad(L)$-orbit.
Let $\hat\O={\cal H}_1\times\ldots\times{\cal H}_r\subset\g_0^r$.
Then by Lemma~5.1 the set of conjugacy classes of subgroups $G(\O)$
for closed $L$-orbits $\O$ in $\hat\O$ is finite.
Let us show that for any $r$-tuple $(x_1,\ldots,x_r)\in\O_1\times\ldots\times\O_r$
there exists an $r$-tuple $(h_1,\ldots,h_r)\in\hat\O$
such that $h_i$ is a homogeneous characteristic of $x_i$
and an $L$-orbit $\Ad(L)(h_1,\ldots,h_r)$ is closed.
Indeed, after simultaneous conjugation of elements $x_i$ by some 
element $g\in L$ we may suppose that any $x_i$ has a homogeneous
characteristic $h_i\in i\k_0$ (in all degrees except at most one
no conjugation is needed because of Moore--Penrose property,
for one degree this is obvious). Then by Lemma~5.2 
an orbit $\Ad(L)(h_1,\ldots,h_r)$ is closed.
Since all functions $\rk\,\Psi_V$ are $L$-invariant,
we may restrict ourselves to the points $x=\sum_i x_i\in\g/\p$
such that $x_i\in\O_i$, any $x_i$ has a homogeneous characteristic
$h_i\in i\k_0$, and a conjugacy class of $\langle h_1,\ldots,h_r\rangle_{alg}$
is fixed. We claim that any function $\rk\,\Psi_V$
is constant along the set of these points. Moreover,
we shall prove that
$$\rk\,\Psi_V(x)=\dim V_R-\dim V_R^{\langle h_1,\ldots,h_r\rangle_{alg}}.
\eqno(*)$$
Indeed, 
$$\rk\,\Psi_V(x)=\dim\sum_i\Im\left(\ad(x_i)|_{V_{R-i}}\right).$$
Clearly $V_R$ is $\ad(h_i)$-invariant and is killed by $\ad(e_i)$, therefore
from the $\ssl_2$-theory we get that
$V_R=\mathop{\oplus}\limits_{k\ge0}V_R^k$, where $\ad(h_i)|_{V_R^k}=k\cdot\Id$.
Moreover, 
$$\Im\left(\ad(x_i)|_{V_{R-i}}\right)=\mathop{\oplus}\limits_{k>0}V_R^k.$$
Let $H$ be a contravariant Hermite form on $V_R$ with respect to the compact form
$\k_0$ of $\g_0$.
Since $h_i\in i\k_0$ and $H$ is a contravariant form we get that
$\mathop{\oplus}\limits_{k>0}V_R^k=(V_R^0)^\perp$.
Therefore,
$$\sum_i\Im\left(\ad(x_i)|_{V_{R-i}}\right)=\left(\cap_i V_R^{h_i}\right)^\perp=
\left(V_R^{\langle h_1,\ldots,h_r\rangle_{alg}}\right)^\perp.$$
The formula $(*)$ follows.
\qed\enddemo

\head \S6. Shortly graded simple real Lie algebras\endhead

Let $\g=\g_{-1}\oplus\g_0\oplus\g_1$ be a real simple shortly graded
Lie algebra. We suppose that $\g$ does not admit a complex structure,
so its complexification $\g^c$ is a simple complex Lie algebra.
There exists a unique element $h\in\g_0$ such that
$ad(c)|_{\g_k}=k\cdot\Id$. Then $C=\R h$ is a center of $\g_0$
and $\g_0=C\oplus\g_0'$, where $\g_0'=[\g_0,\g_0]$.
We fix a maximal compact subalgebra $\k_0\subset\g_0'$ and a Cartan decomposition
$\g_0'=\k_0\oplus\p_0'$. Let $\p_0=\p_0'\oplus C$.
Then we may call  
$\g_0=\k_0\oplus\p_0$ a Cartan decomposition of $\g_0$.

\definition{Definition}
Let $e\in\g_1$. An element $f\in\g_{-1}$ is called a Moore--Penrose
inverse of $e$ if there exists an element $h\in\p_0$ such that
$\langle e,h,f\rangle$ is an $\ssl_2$-triple in $\g$.
\enddefinition

The following is an easy consequence of Theorem 1:

\proclaim{Proposition}
For any $e\in\g_1$ its MP-inverse exists and is unique.
\endproclaim

\demo{Proof}
The complexification $\g^c$ of $\g$ is a shortly graded complex simple
Lie algebra
$$\g^c=\g_{-1}^c\oplus\g_0^c\oplus\g_1^c,\ \text{where}\ \g_k^c=\g_k\otimes\C.$$
Then $\hat\k_0=\k_0\oplus i\p$ is a compact real form of $\g_0^c$.
By Theorem~1 there exists a unique $\ssl_2$ triple $\langle e,h,f\rangle$
in $\g^c$ such that $h\in i\hat\k_0=i\k_0\oplus\p$.
It suffices to show that, in fact, $h\in\p$.
Indeed, $\langle e,\overline h,\overline f\rangle$ is an $\ssl_2$-triple 
such that $\overline h\in i\k_0\oplus\p$, where
bar denotes the complex conjugation in $\g^c$ with respect to $\g$.
By Theorem~1 it follows that $h=\overline h$.
Therefore, $h\in\p$.
\qed\enddemo

The list of shortly graded simple real Lie algebras
is well-known. It can be easily obtained from \cite{D},
where it is shown that there exists a bijection of 
$\Z$-graded real simple Lie algebras and weighted Satake diagrams
with certain natural restrictions.
In the rest part of this section we describe
arising Moore--Penrose inverses.
The proofs are quite straightforward, so they are omitted.
We do this job only for classical real Lie algebras.
For two real forms of $E_6$ and two real forms of $E_7$
that admit short gradings the answer is quite similar to one
obtained in \S3. We only need to change split complex Cayley numbers
to either split real Cayley numbers $\Ca(\R)$
or the division algebra of octonions $\OO$.

\subhead Real and quaternionic linear maps\endsubhead
Suppose that $U=\R^n$ and $V=\R^m$ (resp.~$U=\H^n$ and $V=\H^m$) 
are real vector spaces (resp.~right quaternionic vector spaces)
equipped with standard Euclidean scalar products
(resp.~with standard Hermite scalar products $\sum_i\overline{q_i'}q_i$,
where bar denotes the standard quaternionic involution).
For any linear map $F:\,U\to V$ its Moore-Penrose
inverse is a linear map $F^+:\,V\to U$ defined as follows.
Let $\Ker F\subset U$ and $\Im F\subset V$ be
the kernel and the image of $F$.
Let  $\Ker^\perp F\subset\R^n$ and $\Im^\perp F\subset \R^m$
be their orthogonal complements with respect 
to the Euclidean scalar products (resp.~the Hermite scalar products,
notice that in this case we shall get right vector subspaces).
Then $F$ defines via restriction a bijective linear map
$\tilde F:\,\Ker^\perp F\to\Im F$.
Then $F^+:\,V\to U$ is a unique linear map such that
$F^+|_{\Im^\perp F}=0$ and 
$F^+|_{\Im F}=\tilde F^{-1}$.
This MP-inverse corresponds to short gradings of $\ssl_{n+m}(\R)$ 
(resp.~$\ssl_{n+m}(\H)$).

\subhead Skew-Hermite matrices (or forms)\endsubhead
This example is an analogue of 
a MP-inverse of skew-symmetric forms from the Introduction.
So for diversity we give a matrix description.
We consider matrices over $\R$, $\C$, or $\H$.
The Moore--Penrose inverse of a skew-Hermite 
matrix $A$ (with respect to the canonical
involution, so in the real case $A$ is just skew-symmetric) 
is a unique skew-Hermite 
matrix $A^+$ such that
$$AA^+A=A,\quad A^+AA^+=A^+,\quad [A,A^+]=0.\eqno(*)$$
This MP-inverse corresponds to the short grading of $\so_{p,p}$ (real case),
$\su_{p,p}$ (complex case), $\sp_{p,p}$ (quaternionic case).

\subhead Hermite matrices (or forms)\endsubhead
This example is an analogue of 
a MP-inverse of symmetric forms from the Introduction.
So again we shall give a matrix description.
We consider matrices over $\R$ or $\H$.
The Moore--Penrose inverse of a Hermite 
matrix $A$ is a unique Hermite 
matrix $A^+$ that satisfies equations $(*)$.
This MP-inverse corresponds to the short grading of $\sp_{2p}(\R)$ (real case),
$\u^*_{2p}(\H)$ (quaternionic case).

\subhead Vectors in a pseudo-Euclidean space\endsubhead
We take a vector space $V=\R^{n+m}$ with a standard Euclidean scalar product
$(\cdot,\cdot)$.
Let $I=\left(\matrix \Id_n&0\cr0&-\Id_m\cr\endmatrix\right)$,
where $\Id_k$ is an identity matrix from $\Mat_{k,k}$.
Let $\{u,v\}=(u,Iv)$ be a pseudo-Euclidean scalar product.
The Moore--Penrose inverse takes any vector $v\in V$ 
to a vector $v^\vee$ defined as follows:
$$v^\vee=\cases
-{\textstyle v\over\textstyle \{v,v\}},&\text{if}\ \{v,v\}\ne0\cr
-{\textstyle Iv\over\textstyle 2(v,v)},&\text{if}\ \{v,v\}=0,\ v\ne0\cr
0,&\text{if}\ v=0.\cr
\endcases$$
This MP-inverse corresponds to the short grading of $\so_{n+1.m+1}$.

\head References\endhead
\widestnumber\key{XXXX}

\ref\key CM
\by S.L. Campbell, C.D. Meyer
\paper Generalized inverses of linear transformations
\jour Pitman, London--San Francisco--Melbourne
\yr 1979
\endref

\ref\key D
\by D. Djokovi\'c
\paper Classification of $\Z$-graded real semi-simple Lie
algebras
\jour J.~of Algebra
\vol 76
\yr 1982
\pages 367--382
\endref

\ref\key HR
\by L. Hille, G. R\"ohrle
\paper On parabolic subgroups of classical groups with 
a finite number of orbits on the unipotent radical
\jour C.R. Acad. Sci. Paris
\vol 325, Serie I
\yr 1997
\pages 465--470
\endref

\ref\key Ja
\by N. Jacobson
\paper Structure and representations of Jordan algebras
\jour AMS Colloquium Publications
\vol 39
\yr 1968
\endref

\ref\key Lo
\by O. Loos
\paper Jordan pairs
\jour Lecture notes in Mathematics, Springer--Verlag
\vol 460
\yr 1975
\endref

\ref\key Mo
\by E.H. Moore
\paper General Analysis -- Part I
\jour Mem. Amer. Phil. Soc.
\vol 1
\yr 1935
\endref

\ref\key MRS
\by I. Muller, H. Rubenthaler, G. Schiffmann
\paper Structure des espaces pr\'ehomog\`enes associ\'es \`a certaines
alg\'ebres de Lie gradu\'ees
\jour Math. Ann.
\vol 274
\yr 1986
\pages 95--123
\endref

\ref\key P
\by V.S. Pyasetskii
\paper Linear Lie groups acting with finitely many orbits
\jour Functsional. Anal. i Prilozh.
\vol 9 (4)
\yr 1975
\pages 85--86 (Russian); English. Transl.: Funct. Anal. Appl. {\bf 9} (1975), 351--353
\endref

\ref\key Pa
\by D.I. Panyushev
\paper Parabolic subgroups with abelian unipotent radical as a testing site for invariant
theory
\jour Canad. J. Math.
\vol 51 (3)
\yr 1999
\pages 616--635
\endref

\ref\key Pa1
\by D.I. Panyushev
\paper Complexity and nilpotent orbits
\jour Manuscr. Math.
\vol 83
\yr 1994
\pages 223--237
\endref

\ref\key Pe
\by R.A. Penrose
\paper A generalized inverse for matrices
\jour Proc. Camb. Phil. Soc.
\vol 51
\yr 1955
\endref

\ref\key PV
\by V.L. Popov, E.B. Vinberg
\paper Invariant Theory
\jour in: Algebraic Geometry IV, Encyclopaedia of Mathematical Sciences,
Springer--Verlag
\vol 55
\yr 1994
\pages 123--278
\endref

\ref\key Ri
\by R.W. Richardson
\paper Finiteness theorems for orbits of algebraic groups
\jour Indag. Math.
\vol 88
\yr 1985
\pages 337--344
\endref

\ref\key Ri1
\by R.W. Richardson
\paper Conjugacy classes of $n$-tuples in Lie algebras and
algebraic groups
\jour Duke Math. J.
\vol 57
\yr 1988
\pages 1--35
\endref

\ref\key RRS
\by R.W. Richardson, G. R\"ohrle, and R. Steinberg
\paper Parabolic subgroups with abelian unipotent radical
\jour Invent. Math.
\vol 110
\yr 1992
\pages 649--671
\endref

\ref\key Te
\by E. Tevelev
\paper Isotropic subspaces of multi--linear forms
\jour to appear in Matematicheskie Zametki (in Russian)
\yr 2000
\endref

\ref\key Vi
\by E.B. Vinberg
\paper The Weyl group of a graded Lie algebra
\jour Izv. Akad. Nauk SSSR Ser. Mat.
\vol 40
\yr 1976
\pages 488--526 (Russian); English. Transl.: 
Math. USSR-Izv. {\bf 10} (1976), 463--495
\endref

\ref\key VO
\by E.B. Vinberg, A.L. Onischik
\paper Seminar on Lie Groups and Algebraic Groups
\jour Berlin: Springer
\yr 1990
\endref

\enddocument